\documentclass[11pt]{article}
\usepackage{amssymb,amsmath,latexsym}
\usepackage{mathrsfs}
\usepackage{color}

\def\wt{\widetilde}

\def\<{\left\langle}
\def\>{\right\rangle}
\def\pf{\noindent{\bf Proof.} }

\def\qed{{\hfill $\Box$\medskip}}
\def\to{\rightarrow}

\def\QQ{\mathcal{Q}}

\def\conv{\mathrm{conv}}
\def\PP{\mathbb{P}}
\def\EE{\mathbb{E}}
\def\LL{\mathcal{L}}

\def\VV{\mathcal{V}}
\def\UU{\mathcal{U}}

\def\cvar{\mathrm{CVaR}}
\def\var{\mathrm{VaR}}
\def\RR{\mathbb{R}}
\def\RRR{\mathcal{R}}
\def\PPP{\mathcal{P}}

\def\rv{\sum\limits_{i=1}^na_iX_i}
\def\rob{\sum\limits_{i=1}^na_iz_i}

\input{epsf.sty}

\newtheorem{thm}{Theorem}[section]
\newtheorem{prop}{Proposition}[section]

\newtheorem{example}{Example}[section]
\newtheorem{lemma}{Lemma}[section]

\numberwithin{equation}{section}

\newcommand{\beq}{\begin{equation}}
\newcommand{\eeq}{\end{equation}}
\newcommand{\bay}{\begin{array}}
\newcommand{\eay}{\end{array}}
\newcommand{\bea}{\begin{eqnarray}}
\newcommand{\eea}{\end{eqnarray}}
\newcommand{\beaa}{\begin{eqnarray*}}
\newcommand{\eeaa}{\end{eqnarray*}}

\newcommand{\reff}[1]{(\ref{#1})}

\newcommand{\reft}[1]{Theorem\ \ref{#1}}

\newcommand{\refp}[1]{Proposition\ \ref{#1}}

\def\al{\alpha}

\def\b{\bf}
\def\bdes{\begin{description}}
\def\edes{\end{description}}
\def\benu{\begin{enumerate}}
\def\eenu{\end{enumerate}}
\def\bitm{\begin{itemize}}
\def\eitm{\end{itemize}}
\def\be{\beta}
\def\ba{\bar}
\def\bu{{\scriptsize \Box\,}}

\def\cl{{\rm cl }\,}
\def\cA{{\cal A}}
\def\cd{\cdot}

\def\cs{\cdots}

\def\de{\delta}

\def\E{\mathbb{E}}

\def\esssup{\mathop\mathrm{ess\h{-}sup}}
\def\essinf{\mathop\mathrm{ess\h{-}inf}}

\def\et{\eta}

\def\ga{\gamma}

\def\h{\hbox}

\def\i{\infty}

\def\la{\lambda}

\def\l{\left}

\def\LL{{\sL}}

\def\LR{\ \Longleftrightarrow\ }

\def\n{\noindent}

\def\Om{\Omega}
\def\om{\omega}
\def\ov{\over}
\def\P{\mathbb{P}}

\def\qed{\hspace*{10pt}\hfill{$\square$}\hfilneg\par}
\def\QQb{{\bar{\cal Q}}}
\def\R{\mathbb{R}}

\def\r{\right}

\def\RRb{{\bar{\cal R}}}

\def\sL{\mathscr{L}}

\def\su{\subseteq}
\def\t{\tilde}
\def\ti{\times}
\def\tz{{\t z}}

\begin{document}

\title{
On the Dual Representation of \\
Coherent Risk Measures
\thanks{
Research is partially supported by grants from Australian Research Council (DP160102819) and National University of Singapore. }}
\author{
Marcus Ang\thanks{
Lee Kong Chian School of Business, Singapore Management
University. Email: marcusang@smu.edu.sg} \ \ \
 Jie Sun\thanks{Department of Mathematics and Statistics, Curtin University, Australia, and School of Business, National University of Singapore. Email: jie.sun@curtin.edu.au}\ \ \ Qiang Yao\thanks{School of Statistics, East China Normal University. Email: qyao@sfs.ecnu.edu.cn}  }
\date{}
\maketitle

\thispagestyle{empty}

 \noindent\underline{\hspace*{6.3in}}
\par

{\vskip 10 true pt \noindent{\small{\bf Abstract.}
A classical result in risk measure theory states that every coherent risk measure has a dual representation as the supremum of certain expected value over a risk envelope. We study this topic in more detail.    The related issues include:   1. Set operations of risk envelopes and how they change the risk measures, 2. The structure  of risk envelopes of popular risk measures, 3. Aversity of risk measures and its impact to risk envelopes,   and 4. A connection between risk measures in stochastic optimization and  uncertainty sets in robust optimization.
}

\vskip 14 true pt \noindent{\small\bf Key words.}
{\small  Coherent risk measures, duality,  optimization, risk envelopes }

\noindent\underline{\hspace*{6.3in}}

\vfil\eject
}
\setcounter{page}{1}
\section{Introduction}

At the core of  stochastic optimization is the problem of minimizing $\E_\P[f(x,\tz)],$ where  $x\in\R^n$ is the decision vector, $\tz$ is a random vector, $f:\R^n\ti\R^m\to(-\i,+\i],$  $\E$ stands for expectation, and $\P$ is the joint probability distribution of $\tz$.
In classical numerical stochastic optimization it is  assumed that the distribution of $\P$ is given, which is restrictive since in practice only partial information on $\P$ is available, say, one only knows $\P\in\cA$, where $\cA$ is defined by certain known statistics of $\tz$.   Therefore we are naturally led to a ``distributionally robust" formulation as follows
\begin{eqnarray*}
     {\rm(DRSO)}~~~~&&\min~\sup_{\P\in\cA}\E_\P(f(x,\tz)):=\RRR(f(x,\tz)).
   \end{eqnarray*}
Observe that for a fixed $x$, $X:=f(x,\tz)$ is a random variable and the property of mapping  $\RRR(X)=\sup_{\P\in\cA}\E_\P(X)$  deserves a careful study. In fact, as pointed by Rockafellar (2007),   it is natural to consider the functional $ \RRR(f(x,\tz))$
as a ``risk measure'' or ``surrogate'' of the random cost function $f(x,\tz).$ This paper aims at  studying a dual representation of the function $\RRR$  and its applications in optimization.

Given a probability space $(\Omega, \Sigma, \PP_0)$, it is well known that $X:~\Omega\rightarrow\RR$ is a \emph{random variable} if it is $\Sigma$-measurable, that is, $\{\omega:~X(\omega)\leq a\}\in\Sigma$ for any $a\in\RR$.  We call $\PP_0$ the \emph{base probability measure}, which is  fixed in our analysis. To simplify our notation, when the expectation with respect to $\P_0$  is concerned, we omit $\PP_0$ and write $\EE_{\PP_0}(X)$ as $\EE(X).$ As usual, for $1\leq p\le \infty$, we use $\LL\,^p(\Omega,\Sigma,\PP_0)$~($\LL\,^p$ for short)  to denote the set of all random variables $X$ satisfying $\EE(|X|^p)<+\infty$. For the convenience of engineering applications,  we restrict ourselves to the space of $X\in\sL^2$ although the main results of this paper could be extended to a larger  space such like $\sL^1$. Therefore, in this paper a risk measure $\RRR$ is a functional from $\LL^2$ to $(-\infty,+\infty]$. It may represent ``the risk of loss'' where $X$ may represent ``the real amount of loss''. Furthermore, if $\RRR(X)$ is finite for any $X\in\LL^2$, then we call $\RRR$ a \emph{finite} risk measure.  A risk measure $\RRR$ is {\em coherent in the basic sense} ({``\em coherent''} for short) if it satisfies the following five axioms (Artzner {\it et al.} 1997, 1999, Rockafellar 2007).
\begin{description}
  \item[(A1)] $\RRR(C)=C$ for all constant $C$,
  \item[(A2)] $\RRR((1-\lambda)X+\lambda X')\le (1-\lambda)\RRR (X)+\lambda\RRR(X')$ for $\lambda\in[0,1]$ (``convexity''),
  \item[(A3)] $\RRR(X)\le\RRR(X')$ if $X\le X'$ almost surely (``monotonicity''),
  \item[(A4)] $\RRR(X)\le0$ when $\|X^k-X\|_2\to 0$ with $\RRR(X^k)\le0$ (``closedness''),
  \item[(A5)] $\RRR(\lambda X)=\lambda\RRR(X) $ for $\lambda>0$ (``positive homogeneity'').
\end{description}
In  early literature on coherency (Artzner $et\ al.$ 1997, 1999),  it was required to have $\RRR(X+C)=\RRR(X)+C$. It can be shown that this follows automatically by ({\b A1}) and ({\b A2}) (Rockafellar {\it  et al.} 2006).

Consider another probability measure $\PP$ on $(\Omega,\Sigma)$,  $\PP$ is said to be \emph{absolutely continuous} with respect to  $\PP_0$~(denoted by $\PP\ll\PP_0$) if $\PP_0(A)=0$ implies $\PP(A)=0$ for any measurable set $A\in\Sigma$. If $\PP\ll\PP_0$, then by probability theory there is a well-defined Radon-Nikodym  derivative $Q=\frac{d\PP}{d\PP_0}$. Such derivatives make up the set
\begin{equation}\label{RN}\PPP:=\left\{Q\in\LL^2:~Q\ge0,~\EE(Q)=1\right\}.\end{equation}    $Q$ is called the ``density'' of $\PP$ because the expectation of a random variable $X$ with respect to $\PP$ is equal to $\E(XQ)$, namely
\begin{equation}\label{expect}
\E_{\PP}(X)=\int_\Om X(\om)d\PP(\om)=\int_\Om X(\om)Q(\om)d\PP_0(\om)=\E(XQ).\end{equation}

 Any nonempty closed convex subset $\QQ$ of $\PPP$ is called  a ``risk envelope''. According to the theory of conjugacy in convex analysis, there is a dual representation for coherent risk measures (Theorem 4(a), Rockafellar 2007), which says that
\begin{quote} {\em $\RRR$ is a coherent measure of risk in the basic sense if and only if there is a risk envelope $\QQ$ (which will be uniquely determined) such that}\end{quote}
\begin{equation}\label{DR}\RRR(X)=\sup\limits_{Q\in\QQ}\E(XQ).\end{equation}

Here and below, we will  regard this result as    ``{\em the dual representation theorem}'' for short.

It follows from \reff{DR} that the  risk envelope  $\QQ$    can be written explicitly as
\begin{equation}\label{e:1.1}
\QQ=\{Q\in\PPP:~\EE(XQ)\leq\RRR(X)~\h{for all}~X\in\LL^2\}.
\end{equation}

Note that the requirement $Q\geq0$ in \reff{RN} is equivalent to Axiom ({\b A3})  and the requirement $\EE(Q)=1$ is equivalent to ({\b A1}), as shown in Rockafellar, Uryasev and Zabarankin  (2006). Furthermore, the setting of $X\in \LL^2$ implies $Q\in\sL^2$. Hence all requirements for $Q$ in \reff{RN} are natural. It should be noted that a primary form of the above representation theorem with a finite set $\Omega$ has existed long before the notion of coherent risk measure, see, e.g., Huber (1981).

Many applications of risk measures are concerned with ``averse risk measures''. A risk measure is  \emph{averse} if it satisfies axioms ({\b A1}), ({\b A2}), ({\b A4}), ({\b A5}) and
\begin{description}
  \item[(A6)] $\RRR(X)>\EE(X)$ for all non-constant $X$.
\end{description}
It would be interesting both in theory and practice to  describe aversity in the context of dual representation of risk measures. We shall discuss this topic in Section 4.

The contributions of this paper can be outlined as follows:
\benu
\item We derive formulae of risk measures when the corresponding risk envelopes involve  set operations such as union, intersection, and convex combination (See Proposition 2.1, Theorem 2.1, and Theorem 2.2, respectively).
\item We present independent proofs in Subsections 3.1-3.5 for the correspondence between several popular risk measures and their risk envelopes.
\item We study sufficient and necessary conditions on the risk envelope that guarantee the aversity of the corresponding risk measure (See  Propositions 4.2-4.5).

\item We indicate a connection between the so-called uncertainty sets in robust optimization and the dual representation of risk measures (See Propositions 5.1-5.2specify and Theorem 5.1 for details).
 \eenu

The paper is organized as follows. In Section 2, we consider the set operations of risk envelopes.  In Sections 3 and 4, we  discuss risk envelopes for several popular risk measures and  risk aversity, respectively. Section 5  addresses the relationship between the risk measures defined through uncertainty sets and the ones defined through risk envelopes. Section 6 concludes this paper.

\section{Set Operations of Risk Envelopes }

 Suppose $\RRR_1,\RRR_2,\cdots,\RRR_n$ is a collection of coherent risk measures on $\LL^2$ with risk envelopes $\QQ_1,\QQ_2,\cdots,\QQ_n$ respectively. Since $\LL^2$ is a \emph{Banach lattice}~(that is, it is a Banach space and $X,Y\in\LL^2$ with $|X|\leq|Y|$ implies $\|X\|_2\leq\|Y\|_2$),  if  $\RRR_i$ is finite, then it is continuous, subdifferentiable on $\LL^2$, and   bounded above in some neighborhood of the origin by Proposition 3.1 of Ruszczynski and Shapiro (2006). It then follows that, by Theorem 10 of Rockafellar (1974), the corresponding $\QQ_i$ is compact in the weak topology of $\LL^2$, that is, $\QQ_i$ is \emph{weakly compact}.

The following result deals with convex combination of the sets $\QQ_1,\QQ_2,\cdots,\QQ_n$. A similar result can be found in Rockafellar and Uryasev (2013).

\begin{prop}\label{421} Let $\la_1,...,\la_n$ be positive numbers satisfying $\la_1+\cdots+\la_n=1$. Then the convex combination $$\RRR:=\la_1\RRR_1+\cs+\la_n\RRR_n$$ is a coherent risk measure with risk envelope $$\ba\QQ=\cl(\la_1\QQ_1+\cs+\la_n\QQ_n),$$ where $\cl$ means the closure of the set. Moreover, if all but perhaps one of the $\RRR_i$'s are finite, then  the risk envelope is simply $$\QQ=\la_1\QQ_1+\cs+\la_n\QQ_n.$$
\end{prop}

\pf Since
$$\sup_{Q\in\ba\QQ}\EE(XQ)=\sup_{Q_i\in\QQ_i,i=1,...,n}\EE\l[X(\la_1Q_1+\cs+\la_nQ_n)\r]=\sum_{i=1}^n\la_i\RRR_i(X)=\RRR(X),$$
the first part of the proposition follows. For the second part, as discussed above, we know that if  $\RRR_i$ is finite, then the corresponding $\QQ_i$ is weakly compact. It is easy to see that $\QQ$ is a nonempty and convex subset of $\PPP$~(as defined in (\ref{RN})). Furthermore, $\QQ$ is weakly closed since all but perhaps one of the $\QQ_i$'s are weakly compact, and the sum of finitely many weakly closed set, if all but perhaps one of which is weakly compact, is a weakly closed set. Then $\QQ$ is closed because closedness coincides with weak closedness for convex sets. Therefore, $\ba{\QQ}=\QQ$ in this case. \qed

Next, define
\beaa
&\widetilde{\RRR}_1(X)
:=\max\limits_{1\leq i\leq n}\RRR_i(X),\quad\widetilde{\RRR}_2(X)
:=\min\limits_{1\leq i\leq n}\RRR_i(X),\h{ and}\\
&\widetilde{\RRR}_3(X):=\cl (\RRR_1\bu\RRR_2\bu\cdots\Box\RRR_n)(X),
\eeaa
where $\cl$ means the closure of the function (Rockafellar and Wets 1997)  and
$$(\RRR_1\bu\RRR_2\bu\cdots\bu\RRR_n)(X):=\inf\{\RRR_1(X_1)+\RRR_2(X_2)+\cdots+\RRR_n(X_n):~X_1+X_2+\cdots+X_n=X\}$$ is the so-called {\em inf-convolution} of the functionals $\RRR_i, i=1,...,n.$
Let us call $\widetilde{\RRR}_1 $ and $\widetilde{\RRR}_2$ the ``max'' and the ``min'' of the risk measures $\RRR_1,\RRR_2,\cdots,\RRR_n$, respectively. Clearly,  $\widetilde{\RRR}_2(X)$ is not  coherent because it may not be convex. We next show that $\widetilde\RRR_1$ and the lower-convexification of $\widetilde\RRR_2$, namely $\widetilde\RRR_3$, are coherent risk measures generated by the risk envelopes $\conv\left(\bigcup\limits_{i=1}^n\QQ_i\right)$ and $\bigcap\limits_{i=1}^n\QQ_i$, respectively, where $\conv(\cd)$ stands for the convex hull. We  begin with the following  lemma about $\widetilde{\RRR}_2$ and $\widetilde{\RRR}_3$.

\begin{lemma}\label{p:1.2}
$\widetilde{\RRR}_3$ is the ``lower-convexification'' of $\widetilde{\RRR}_2$ in the sense that

(1) $\widetilde{\RRR}_3(X)\leq\widetilde{\RRR}_2(X)$ for all $X$.

(2) Let $\RRR(X)$ be any coherent risk measure  satisfying
$\RRR(X)\leq\widetilde{\RRR}_2(X)$ for all $X$. Then $\RRR(X)\leq\widetilde{\RRR}_3(X)$ for all $X$.

\end{lemma}
\pf (1) By the definition of $\widetilde{\RRR}_3$, we have for any $1\leq i\leq n$ and for all $X$,
$$\widetilde{\RRR}_3(X)\le\cl\big[\RRR_1(0)+\cdots+\RRR_{i-1}(0)+\RRR_i(X)+\RRR_{i+1}(0)+\cdots+\RRR_n(0)\big]=\RRR_i(X).$$
Then  $\widetilde{\RRR}_3(X)\leq\min\limits_{1\leq i\leq n}\RRR_i(X)=\widetilde{\RRR}_2(X)$  as desired.

(2) Since $\RRR(X)\leq\widetilde{\RRR}_2(X)$ for all $X$, we have $\RRR(X)\leq\RRR_i(X)$ for any $1\leq i\leq n$ and for all $X$. Furthermore, by the convexity of $\RRR$, we have for any $X_1,X_2,\cdots,X_n$ such that $X_1+X_2+\cdots+X_n=X$,
$$\RRR(X)\leq\RRR(X_1)+\RRR(X_2)+\cdots+\RRR(X_n)\leq\RRR_1(X_1)+\RRR_2(X_2)+\cdots+\RRR_n(X_n).$$
Taking closure of  infimum on the right hand side, by the definition of $\widetilde{\RRR}_3$ together with the continuity of $\RRR_1,\cdots,\RRR_n$, we get $\RRR(X)\leq\widetilde{\RRR}_3(X)$ for all $X$, as desired.\qed

 The main results of this section are the following two theorems. A finite-dimensional version of them appeared in Theorem 3.3.3 of Hiriart-Urruty and Lemarach\'el (1993).  Here, we present a proof for the $\sL^2$ version.
\begin{thm}\label{t:1.1}If $\RRR_1,\cdots,\RRR_n$ are finite, then $\wt\RRR_1(\cdot)$ is a coherent risk measure with risk envelope
$\widetilde{\QQ}_1=\conv\left(\bigcup\limits_{i=1}^n\QQ_i\right)$.
\end{thm}
\pf We first claim that $\conv\left(\bigcup\limits_{i=1}^n\QQ_i\right)$ is closed and convex. The convexity is trivial. For closedness, since $\QQ_1,\cdots,\QQ_n$ are all weakly compact, we have that $\conv\left(\bigcup\limits_{i=1}^n\QQ_i\right)$ is weakly compact because the union of any finite collection of weakly compact sets is again weakly compact, and its convex hull is therefore weakly compact. Furthermore, $\conv\left(\bigcup\limits_{i=1}^n\QQ_i\right)$ is closed because weak compactness implies weak closedness, and weak closedness coincides with closedness for convex sets. Next, for any $X\in\LL^2$, we have
$$\widetilde{\RRR}_1(X)=\max\limits_{1\leq i\leq n}\RRR_i(X)=\max\limits_{1\leq i\leq n}\left(\sup\limits_{Q\in\QQ_i}\EE(XQ)\right)=\sup\limits_{Q\in\bigcup\limits_{i=1}^n\QQ_i}\EE(XQ)=\sup\limits_{Q\in\conv\left(\bigcup\limits_{i=1}^n\QQ_i\right)}\EE(XQ).$$
Hence by the dual representation theorem, $\widetilde{\RRR}_1$ is a coherent risk measure and its risk envelope is  $\widetilde{\QQ}_1=\conv\left(\bigcup\limits_{i=1}^n\QQ_i\right)$, as desired.\qed

\begin{thm}\label{t:1.2}
$\widetilde{\RRR}_3(\cdot)$ is a coherent risk measure with risk envelope $\bigcap\limits_{i=1}^n\QQ_i$ if and only if $\bigcap\limits_{i=1}^n\QQ_i\neq\emptyset$.
\end{thm}

\pf For the ``if'' part, we first verify that $\widetilde{\RRR}_3(\cdot)$ is a coherent risk measure. By the closure of inf-convolution formula of $\widetilde{\RRR}_3$, the convexity ({\b A2}) and closedness ({\b A4}) hold. For  positive homogeneity ({\b A5}), one has
\beaa
\wt\RRR_3(\la X)&=&\cl\inf_{X_2,...,X_n}\l\{\RRR_1(\la X-X_2-...-X_n)+\RRR(X_2)+\cs+\RRR_n(X_n)\r\}\\
&=&\cl\inf_{Y_2,...,Y_n}\l\{\RRR_1(\la X-\la Y_2-...-\la Y_n)+\RRR(\la Y_2)+\cs+\RRR_n(\la Y_n)\r\}\\
&=&\la\wt\RRR_3(X).
\eeaa
Axiom ({\b A1}) is true because
\beq\label{c1}\wt\RRR_3(C)\le \RRR_1(C)+\RRR_2(0)+\cs+\RRR_n(0)= C\h{ and similarly, }\wt\RRR_3(-C)\le-C.\eeq
Then by convexity and positive homogeneity
\beq\label{c2}0 =\wt \RRR_3(0)\le\wt\RRR_3(C) +\wt\RRR_3(-C)\le\wt\RRR_3(C)- C\ \LR \ \wt\RRR_3(C)\ge C.\eeq
Thus, ({\b A1}) follows. Finally, let $X\le Y$ almost surely. Then
\beaa
\wt\RRR_3(X)&=&\cl\inf_{X_2,...,X_n}\l\{\RRR_1(X-X_2-...-X_n)+\RRR(X_2)+\cs+\RRR_n(X_n)\r\}\\
&\le&\cl\inf_{X_2,...,X_n}\l\{\RRR_1(Y-X_2-...-X_n)+\RRR(X_2)+\cs+\RRR_n(X_n)\r\}\\
&=&\wt\RRR_3(Y),
\eeaa
hence monotonicity ({\b A3}) holds. Therefore, $\wt\RRR_3(X)$ is a coherent risk measure. Let $\wt\QQ_3$ be its risk envelope. Since $\wt\RRR_3(X)\le\RRR_i(X),$ by \reff{e:1.1}, $\wt\QQ_3\su\QQ_i$ for $1\le i\le n.$ Thus, $\wt\QQ_3\su \bigcap\limits_{i=1}^n\QQ_i$. Conversely, suppose $\widetilde{\RRR}$ is the risk measure with envelope $\bigcap\limits_{i=1}^n\QQ_i$. Since $\widetilde{\RRR}$ is convex, positive homogeneous, and $\widetilde{\RRR}(X)\leq\widetilde{\RRR}_2(X)$ for all $X$, by Lemma \ref{p:1.2} we get $\widetilde{\RRR}(X)\leq\widetilde{\RRR}_3(X)$ for all $X$. Using (\ref{e:1.1}) again, we can get $\bigcap\limits_{i=1}^n\QQ_i\subseteq\widetilde{\QQ}_3$. Thus, we have $\widetilde{\QQ}_3=\bigcap\limits_{i=1}^n\QQ_i$.

We next prove the ``only if'' part. If $\widetilde{\RRR}_3(\cdot)$ is a coherent risk measure, then it has a nonempty risk envelope $\widetilde{\QQ}_3$, which is an implication of Axiom ({\b A1}) and the dual representation theorem. Using the same argument from the last paragraph, we can get $\widetilde{\QQ}_3\subseteq\bigcap\limits_{i=1}^n\QQ_i$. Therefore, $\bigcap\limits_{i=1}^n\QQ_i\neq\emptyset$. \qed
\bigskip

Note that \reft{t:1.2} does not require the $\RRR_i$s to be finite.

Set operations of risk envelopes may be used to create new risk measures that are more conservative (say, by union) or more aggressive (say, by intersection) in applications. Chen $et\ al.$ (2010) used intersections of five uncertainty sets to create new uncertainty sets in robust optimization and here we have shown the same principle applies to risk envelopes.

\section{Popular risk measures and their risk envelopes}\label{section:examples_coherent}

Besides set operations, one can create various different coherent risk measures by adding additional functional constraints to the risk envelope $\PPP$ in \reff{RN}. In this section we study 1) risk measure from expectation, 2) risk measure from worst case analysis, 3) risk measure
from subdividing the future,  4) risk measures from the conditional value at risk and optimized certainty equivalence, and 5)
risk measure from mean-deviation.  Most of the results in this section have been stated in Rockafellar (2007) without proofs. In fact their proofs are scattered in the literature via different approaches.  Here we provide independent proofs based on the unified view of dual representation of risk measures. Our approach is to directly specify the risk envelope $\QQ$ for each of the above cases and to verify the relationship $\RRR(X)=\sup\limits_{Q\in\QQ}\E(XQ).$ The coherency of $\RRR$ then follows from the dual representation theorem.

\subsection{Risk envelope for  expectation}\label{expectation}

Here $\QQ=\{Q\in\sL^2:Q\equiv1\}.$ Then $\EE(X)=\sup_{Q\in\QQ}\EE(XQ).$

\subsection{Risk envelope for the worst case}\label{esssup}

Here the risk envelope is $\QQ=\PPP$ and by ``the worst case'' we mean the   ``essential supremum'' function of $X$, that is,
\begin{equation}\label{e:esssupdef}
\esssup(X):=\inf\{a:~\PP_0(X>a)=0\}.
\end{equation}
Note that $\sup\limits_{Q\in\PPP}\EE(XQ)\leq\esssup(X)$ for any $X\in\LL^2$, and therefore $\PPP\subseteq\QQ$. Hence $\QQ=\PPP$.

It is possible that $\esssup(X)=\i$ for some $X$, which could happen  if $X$ does not have a finite essential supremum. Thus, $\esssup(\cd)$ is not a finite risk measure.

\subsection{The risk measure from subdividing the future }\label{dividing}
In Rockafellar (2007) the following risk measure is discussed. Let $\Omega$ be partitioned into subsets $\Omega_1,\cdots,\Omega_r, r\ge2,$ having positive probability $\PP_0(\Om_k)=\la_k$ with $ \la_1+\cs+\la_r=1.$ For $k=1,\cdots,r$, let
$$\RRR_k(X):=\esssup\limits_{\omega\in\Omega_k}X(\omega):=\inf\{a:~\PP_0(\{X>a\}\cap\Omega_k)=0\}.$$
Then
\begin{equation}\label{dividing_formula}
\RRR:=\lambda_1\RRR_1+\cdots+\lambda_r\RRR_r
\end{equation}
is  a coherent risk measure,  called the risk measure from subdividing the future, whose risk envelope is
\begin{equation}\label{dividing_envelope}
\QQ:=\lambda_1\QQ_1+\cdots+\lambda_r\QQ_r\quad\h{with }\QQ_k:=\{Q\in\PPP:~\EE(Q\textbf{1}_{\Omega_k})=1\}.
\end{equation}

To prove this by Proposition \ref{421}, we only need to prove that $\QQ$ is closed.  Suppose $Q_n\in\lambda_1\QQ_1+\cdots+\lambda_r\QQ_r$ for $n=1,2,\cdots$ and $\|Q_n-Q\|_2\rightarrow0$ as $n\rightarrow\infty$. Then by (3.3), for $n=1,2,\cdots$ we have $\EE(Q_n\textbf{1}_{\Omega_k})=\lambda_k$ for $k=1,2,\cdots,r$. Note that for $k=1,2,\cdots,r$, $$|\EE(Q_n\textbf{1}_{\Omega_k})-\EE(Q\textbf{1}_{\Omega_k})|\leq\|Q_n-Q\|_2\cdot[\PP_0(\Omega_k)]^{\frac{1}{2}}\rightarrow0$$ as $n\rightarrow\infty$. Thus, $\EE(Q\textbf{1}_{\Omega_k})=\lambda_k$ for $k=1,2,\cdots,r$, and therefore $Q\in\lambda_1\QQ_1+\cdots+\lambda_r\QQ_r$. This implies $\lambda_1\QQ_1+\cdots+\lambda_r\QQ_r$ is closed in $\LL^2$.
\qed

\subsection{The conditional value at risk (CVaR) and the optimized certainty equivalence (OCE)}\label{oce}
An important coherent risk measure is the conditional value at risk, popularized by Rockafellar and Uryasev (2000), with the formula
\beq\label{32}{\rm CVaR}_\alpha(X)=\min_{\beta\in\RR}\l\{  \beta+{1\over 1-\alpha}\EE(X-\beta)_+\r\},\eeq where $(t)_+=\max(t,0)$.
We next prove that the risk envelope of CVaR is
$$\QQ_\alpha:=\left\{Q\in\sL^2:~\E(Q)=1,0\le Q\le{1\over1-\alpha}\right\}.$$

 For any $Q\in\QQ_\alpha$ and $\beta\in\RR,$ we have
\begin{eqnarray*}
\EE(XQ)&=&\EE\left[(X-\beta)Q\right]+\beta\EE(Q)\\
&\le&\beta+\EE[Q(X-\beta)_+]\le\beta+{1\over 1-\alpha}\EE(X-\beta)_+.
\end{eqnarray*}
Taking supremum on the left hand side over $Q\in\QQ_\alpha$ and infimum on the right hand side over all $\beta\in\RR,$ we get
\begin{equation}\label{01}\sup_{Q\in\QQ_\al}\EE(XQ)\le \min_\beta\l\{\beta+{1\over 1-\alpha}\EE(X-\beta)_+\r\}.
\end{equation}
On the other hand, noting that the ``value-at-risk'' (VaR) is defined as
$$\var_\al(X):=\inf\l\{\nu\in\R:\P(X>\nu)<1-\al\r\},$$ we have
$$\PP_0(X>\hbox{VaR}_\alpha(X))\le1-\alpha\le\PP_0(X\ge \hbox{VaR}_\alpha(X)).$$
Thus, there exists $\lambda\in[0,1]$ such that $$1-\alpha=\lambda\cdot\PP_0(X>\var_\alpha(X))+(1-\lambda)\cdot\PP_0(X\geq\var_\alpha(X)).$$
Set $$Q_0=\frac{1}{1-\alpha}\cdot[\lambda\cdot\textbf{1}_{\{X>\var_\alpha(X)\}}+(1-\lambda)\cdot\textbf{1}_{\{X\geq\var_\alpha(X)\}}].$$
Note that $0\leq Q_0\leq\frac{1}{1-\alpha}$ and $\EE(Q_0)=1$. Thus $Q_0\in\QQ_\alpha$ and
\begin{eqnarray*}
\sup_{Q\in\QQ_\alpha}\EE(XQ)&\ge&\EE(XQ_0)\\
&=&\EE[(X-\var_\alpha(X))\cdot Q_0]+\var_\alpha(X)\cdot\EE(Q_0)\\
&=&{\rm VaR}_\alpha(X)+{1\over1-\alpha}\cdot\EE(X-{\rm VaR}_\alpha(X))_+\\
&\ge&\min_{\beta\in\RR}\left\{  \beta+{1\over 1-\alpha}\EE(X-\beta)_+\right\}.
\end{eqnarray*}
Combine (\ref{01}) and the above we obtain that
$$\cvar_\al(X)=\sup_{Q\in\QQ_\alpha}\EE(XQ).$$
As a by-product of the proof,  we see that  the minimum in (\ref{32}) is attained at $\beta=\var_\alpha(X)$, that is,
$$\cvar_\alpha(X)=\var_\alpha(X)+\frac{1}{1-\alpha}\cdot\EE\l(X-\var_\alpha(X)\r)_+.$$

Ben-Tal and Teboulle (2007) proved that the negative of their OCE function
$${\rm OCE}_u(X)=\sup_\et\{\et+\E[u(X-\et)]\},$$
where $u$ is a piecewise linear utility function, is a coherent risk measure that includes CVaR as a special case. Since $X$ is a risk rather than an income in our context and we are considering risk rather than utility, we define
\beq\label{001}S_r(X):=-{\rm OCE}_u(-X)=\inf_\et\{-\et+\E[-u(-X-\et)]\}=\inf_\be\{\be+\E[r(X-\be)]\},\eeq
where $r(X)=-u(-X)$ and we can similarly show that if
$$r(X)=\ga_1[X]_+-\ga_2[-X]_+ \h{ with }0\le\ga_2<1<\ga_1,$$ then $S_r(X)$ is a coherent risk measure with  risk envelope $\ga_2\le Q\le \ga_1.$ i.e.,
\beq\label{33}S_r(X)=\sup_{Q\in\QQ_{\ga_1,\ga_2}}\EE(XQ),\hbox{ where }\QQ_{\ga_1,\ga_2}:=\left\{Q\in\PPP:~\ga_2\le Q\le\ga_1\right\}.\eeq

It is interesting to observe that   OCE can be representable by CVaR, namely
$$S_r(X)=\ga_2\E(X)+\cvar_\al(X),\h{ where }\al=1-(\ga_1-\ga_2)^{-1}.$$
This formula can be obtained by using \refp{421} and the fact
$$Q_{\ga_1,\ga_2}=\ga_2\{1\}+Q_\al.$$

\subsection{The mean-deviation }\label{meandeviation}

Fix $0\leq\lambda\leq1$. Define the mean-deviation risk measure as $$\RRR(X)=\EE X+\lambda\cdot\|(X-\EE X)_+\|_2$$ for all $X\in\sL^2$, where $\|\cdot\|_2$ denotes the $\sL^2$-norm, that is, $\|X\|_2:=\left[\EE(X^2)\right]^{\frac{1}{2}}.$

Similar to \reff{e:esssupdef}, we define
\begin{equation}\label{e:essinfdef}
\essinf(X):=\sup\{a:~\PP_0(X<a)=0\}.
\end{equation}
We claim that the risk envelope of $\RRR$ is $$\QQ=\left\{0\leq Q\in\LL^2:~\EE(Q)=1,~\|Q-\essinf Q\|_2\leq\lambda\right\}.$$ In fact, on one hand, for any $X\in\LL^2$ and $Q\in\QQ$, we have
\beaa
\EE(XQ)&=&\EE[(X-\EE X)(Q-\essinf Q)]+\EE X\leq\EE X+\EE[(X-\EE X)_+(Q-\essinf Q)]\\
&\leq&\EE X+\|(X-\EE X)_+\|_2\cdot\|Q-\essinf Q\|_2\leq\EE X+\lambda\cdot\|(X-\EE X)_+\|_2
\eeaa
by Cauchy-Schwartz inequality. Hence we get
\begin{equation}\label{e:3.5}
\sup\limits_{Q\in\QQ}\EE(XQ)\leq\EE X+\lambda\cdot\|(X-\EE X)_+\|_2
\end{equation}
for any $X\in\LL^2$. On the other hand, set
$$Q_0:=1+\frac{\lambda\cdot\left[(X-\EE X)_+-\EE(X-\EE X)_+\right]}{\|(X-\EE X)_+\|_2}.$$
Since $0\leq\lambda\leq1$, we have
$$\essinf Q_0=1-\frac{\lambda\cdot\EE(X-\EE X)_+}{\|(X-\EE X)_+\|_2}\geq1-\frac{\EE(X-\EE X)_+}{\|(X-\EE X)_+\|_2}\geq0.$$
Thus, $0\leq Q_0\in\LL^2$, $\EE Q_0=1$ and
$$\|Q_0-\essinf Q_0\|_2=\frac{\left\|\lambda\cdot(X-\EE X)_+\right\|_2}{\|(X-\EE X)_+\|_2}=\lambda,$$
that is, $Q_0\in\QQ$. Then for any $X\in\LL^2$,
\bea\label{e:3.6}
\sup\limits_{Q\in\QQ}\EE(XQ)&\geq&\EE(XQ_0)=\EE X+\frac{\lambda\cdot\EE\left[(X-\EE X)_+\cdot(X-\EE X)\right]}{\|(X-\EE X)_+\|_2}\nonumber\\
&=&\EE X+\frac{\lambda\cdot\|(X-\EE X)_+\|_2^2}{\|(X-\EE X)_+\|_2}=\EE X+\lambda\cdot\|(X-\EE X)_+\|_2.
\eea
(\ref{e:3.5}) and (\ref{e:3.6}) together imply
$$\sup\limits_{Q\in\QQ}\EE(XQ)=\EE X+\lambda\cdot\|(X-\EE X)_+\|_2.$$
We can check that $\QQ$ is nonempty, convex and closed in $\LL^2$. Therefore, it is the risk envelope for the mean-deviation risk measure.

It should be noted that $\la\le1$ is necessary for coherency as shown by the following example. Consider $$\RRR(X)=\EE X+\lambda\cdot\|(X-\EE X)_+\|_2,$$ where $X$ is a discrete random variable with distribution $$\PP(X=-1)=p,~~~~~~~\PP(X=0)=1-p,$$ where $0<p<1$. Then $\EE X=-p$, so $$\PP((X-\EE X)_+=0)=p,~~~~~~~~\PP((X-\EE X)_+=p)=1-p,$$ and therefore, $\RRR(X)=-p+\lambda p\sqrt{1-p}=p(\lambda\sqrt{1-p}-1)$. If $\lambda>1$, we can take $p>0$ sufficiently small to get $\RRR(X)>0$. However, since we have $X\leq0$ almost surely, this contradicts monotonicity.

\section{Discussion on Aversity}

In this section, we study the effect of aversity on risk measures.
Suppose $\RRR$ is a functional from $\LL^2$ to $(-\infty,+\infty]$. Recall that an \emph{averse} risk measure is defined by axioms ({\b A1}), ({\b A2}), ({\b A4}), ({\b A5}) and
\begin{description}
  \item[(A6)] $\RRR(X)>\EE(X)$ for all non-constant $X$.
\end{description}

 We are interested in the risk measures which are both coherent and averse. Next we  develop the conditions of risk envelopes under which a coherent risk measure is averse. We use the notion ``$A\subset B$'' to denote that $A$ is a proper subset of $B$, that is, $A\subseteq B$ but $A\neq B$. The following necessary condition is trivial.
\begin{prop}\label{p:averse_ness}
Suppose $\RRR$ is a coherent risk measure on $\LL^2$ with risk envelope $\QQ$. If $\RRR$ is averse, then $\{{\b 1}\}\subset\QQ$.
\end{prop}
On the other hand, a sufficient condition is stated in the following proposition.
\begin{prop}\label{p:averse_suff}
Suppose $\RRR$ is a coherent risk measure  with risk envelope $\QQ$. If ${\b 1}$ is a relative interior point of $\QQ$~(relative to $\PPP$), then $\RRR$ is averse.
\end{prop}
\pf Since ${\b 1}$ is a relative interior point of $\QQ$~(relative to $\PPP$), there exists $\delta\in(0,1)$ such that
\begin{equation}\label{e:interior}
\{Q\in\PPP:~\|Q-{\b 1}\|_2<\delta\}\subseteq\QQ.
\end{equation}
If $X$ is not a constant almost surely, then there exists $b\in\RR$ such that $$\PP_0(X\geq b)=p\in(0,1),~~~~~~\PP_0(X<b)=1-p\in(0,1).$$ Set
$$Q_0:=\left\{\begin{array}{ll}1+(1-p)\delta~~~~~~~\h{if}~X\geq b,\\~~~1-p\delta~~~~~~~~~~~~\h{if}~X<b. \end{array}\right.$$ Then we have
$$Q_0\geq0,~~~\EE(Q_0)=1,~~~\|Q_0-{\b 1}\|_2<\delta.$$ By (\ref{e:interior}), we can get that $Q_0\in\QQ$. Thus,
\begin{equation}\label{ineq1}
\EE(XQ_0)\leq\sup\limits_{Q\in\QQ}\EE(XQ)=\RRR(X).
\end{equation}
Furthermore, we have
\bea\label{ineq2}
\EE(XQ_0)-\EE(X)&=&(1-p)\delta\cdot\EE(X\textbf{1}_{\{X\geq b\}})-p\delta\cdot\EE(X\textbf{1}_{\{X<b\}})\nonumber\\
&>&(1-p)\delta b\cdot\PP_0(X\geq b)-p\delta b\cdot\PP_0(X<b)=0.
\eea
(\ref{ineq1}) and (\ref{ineq2}) together imply that $\RRR(X)>\EE(X)$ for all non-constant $X$. Therefore, $\RRR$ is averse.\qed

\bigskip

From Propositions \ref{p:averse_ness} and \ref{p:averse_suff}, we can get the following:
\begin{equation}\label{averse}
{\b 1}~\h{is a relative interior point of}~\QQ~(\h{relative to}~\PPP)\Longrightarrow\RRR~\h{is averse}\Longrightarrow\{{\b 1}\}\subset\QQ.
\end{equation}

Generally, the converse of (\ref{averse}) may not be true, which can be seen from the following two examples.

\begin{example}\label{exam:averse1}
Suppose $\Omega=[0,1]$, $\Sigma$ is the Borel sigma algebra on $[0,1]$, and $\PP_0$ is the Lebesgue measure. In this case $$\{{\b 1}\}:=\{{\t Q}_1(\om)\equiv1\}.$$Consider $\RRR=\cvar_{0.5}$. By Rockafellar (2007), $\RRR$ is a coherent and averse risk measure with risk envelope $\QQ=\{Q\in\LL^2:~0\leq Q\leq2,~\EE(Q)=1\}$. However, ${\b 1}$ is not a interior point of $\QQ$. In fact, for any $\delta\in(0,1)$, the random variable $\widetilde{Q}_\de$ defined as
$$\widetilde{Q}_\de(\omega)=\left\{\begin{array}{ll}~~~3~~~~~~~~\omega\in\left[0,\frac{\delta^2}{16+\delta^2}\right],\\1-\frac{\delta^2}{8}~~~~~\omega\in\left(\frac{\delta^2}{16+\delta^2},1\right] \end{array}\r.$$ is arbitrarily close to $\widetilde{Q}_1(\om)$, but $\widetilde{Q}_\de\not\in\QQ$. Therefore, ${\b 1}$ is not a relative interior point of $\QQ$. Hence the converse of the first ``$\Longrightarrow$''in \reff{averse}  may not be true.
\end{example}

\begin{example}\label{exam:averse2}
Suppose $\Omega=\{\omega_1,\omega_2,\omega_3\}$ and $\PP_0(\{\omega_1\})=\PP_0(\{\omega_2\})=\PP_0(\{\omega_3\})=1/3$. Let
$$Q_0:~~Q_0(\omega_1)=\frac{3}{4},~Q_0(\omega_2)=\frac{3}{2},~Q_0(\omega_3)=\frac{3}{4}.$$ Then $Q_0\in\PPP$ and in this case
$${\b 1}:=Q_1:~~Q_1(\omega_1)=1,~Q_1(\omega_2)=1,~Q_1(\omega_3)=1.$$
 Take $\QQ:=\conv\{Q_1,Q_0\}$, then $\{{\b 1}\}\subset \QQ$. However, for the non-constant random variable $$X:~X(\omega_1)=-1,~X(\omega_2)=0,~X(\omega_3)=1,$$ one has
$$\RRR(X)=\sup_{Q\in\QQ}\EE(XQ)=\max\{\EE(XQ_1),\EE(XQ_0)\}=0=\EE(X).$$ Therefore, $\RRR$ is not averse.
\end{example}

From Example \ref{exam:averse2} we can see that the converse of the second ``$\Longrightarrow$'' in \reff{averse} may not hold even when $\Omega$ is finite. However, the converse of the first ``$\Longrightarrow$'' always holds when $\Omega$ is finite, see the following proposition.
\begin{prop}\label{p:averse_disc}
If $\Omega$ is finite and $\RRR$ is a coherent risk measure  with risk envelope $\QQ$, then $\RRR$ is averse if and only if ${\b 1}$ is a relative interior point of $\QQ$.
\end{prop}
\pf By Proposition \ref{p:averse_suff}, we only need to prove one direction, that is, aversity implies that ${\b 1}$ is a relative interior point. Suppose $\Omega=\{\omega_1,\cdots,\omega_n\}$ and $\PP_0(\{\omega_i\})=p_i>0$ for $i=1,2,\cdots,n$. In this case, $$\PPP=\left\{(q_1,\cdots,q_n):~q_1,\cdots,q_n\geq0,~\sum\limits_{i=1}^nq_ip_i=1\right\},$$ and the risk envelope of $\RRR$ is certain nonempty closed convex $\QQ\subseteq\PPP$, that is, $$\RRR(X)=\max\limits_{(q_1,\cdots,q_n)\in\QQ}\{x_1q_1p_1+\cdots+x_nq_np_n\}$$ for $X=(x_1,\cdots,x_n)\in\RR^n$. Here, $x_i=X(\omega_i)$ for $i=1,2,\cdots,n$. Moreover, since $\RRR$ is averse, we have
\begin{equation}\label{e:averse}
\max\limits_{(q_1,\cdots,q_n)\in\QQ}\{x_1q_1p_1+\cdots+x_nq_np_n\}>x_1p_1+\cdots+x_np_n
\end{equation}
whenever $X(\om_i)$ is not a constant. Note that the affine hull of $\PPP$ is a hyperplane of dimension $n-1$ with a normal vector $(p_1,...,p_n)$. Let the apostrophe of a vector represent its transpose.  Therefore, to prove that $(1,\cdots,1)$ is an interior point of $\QQ$ relative to $\PPP$, we only need to prove that
 \beq\label{46}\max\limits_{(q_1,\cdots,q_n)\in\QQ}(y_1,\cdots,y_n)[(q_1,\cdots,q_n)-(1,\cdots,1)]'>0\eeq for any $(y_1,\cdots,y_n)$ that is not a normal vector of the affine hull of $\PPP$.  In other words, we show that \reff{46} holds for any $(y_1,\cdots,y_n)$ that is not a  multiple of $(p_1,\cdots,p_n)$.

To prove  \reff{46}, noting that
if  $\frac{y_1}{p_1},\cdots,\frac{y_n}{p_n}$ are not the same, then setting $x_i={y_i\ov p_i}$ in (\ref{e:averse}),  we have
\beaa\max\limits_{(q_1,\cdots,q_n)\in\QQ}\{y_1q_1+\cdots+y_nq_n\}
&=&\max\limits_{(q_1,\cdots,q_n)\in\QQ}\left\{\frac{y_1}{p_1}\cdot q_1p_1+\cdots+\frac{y_n}{p_n}\cdot q_np_n\right\}\\
&=& \max\limits_{(q_1,\cdots,q_n)\in\QQ}\{x_1q_1p_1+\cdots+x_nq_np_n\}\\
&>&x_1p_1+\cdots+x_np_n\\
&=&y_1+\cdots+y_n.\eeaa
Therefore \reff{46} is true, implying that $(1,1,\cdots,1)$ is an interior point of $\QQ$ relative to $\PPP.$\qed

We next analyze the examples in Section \ref{section:examples_coherent}. Obviously, the expectation measure $\EE(\cdot)$ in subsection \ref{expectation} is not averse. We call a risk measure $\RRR$ ``law-invariant'' if $\RRR(X)=\RRR(Y)$ whenever $X$ and $Y$ have the same distribution under $\PP_0$. F\"ollmer and Schied (2002)  proved that if $\RRR$ is a coherent, law-invariant risk measure in $\LL^\infty$~(not $\LL^2$) other than $\EE(\cdot)$, then  $\RRR$ is averse. Therefore,  the examples in subsections  \ref{esssup}, \ref{oce} and \ref{meandeviation} are all averse.  However, since we are considering the $\LL^2$ case, we cannot use the result in F\"ollmer and Schied (2002)  directly. We also  noted that the result in $\LL^2$ space has appeared in Rockafellar and Uryasev (2013) without proof. For completeness, we give a direct proof in the next proposition.

\begin{prop}\label{p:averse_examples}
The worst-case, CVaR, OCE and mean-deviation, as risk measures, are all averse.
\end{prop}
\pf The proof is trivial for $\esssup(\cdot)$, since the expectation of any random variable is no larger than its essential supremum, and they are equal if and only if the random variable is a constant almost surely.

For the mean deviation measure, obviously, we have $\EE X+\lambda\cdot\|(X-\EE X)_+\|_2\geq\EE X$ for any $X\in\LL^2$, in which the equality holds if and only if $X\leq\EE X$ almost surely, which implies $X=\EE X$~(i.e. $X$ is a constant) almost surely. Therefore, the mean deviation measure is averse.

For the OCE measure, since $1\in\QQ_{\ga_1,\ga_2}$, we have $S_r(X)\geq\EE(X)$ by \refp{p:averse_suff}. Next, if
$$\EE(X)=S_r(X)=\min_{\beta\in\RR}\big\{  \beta+\EE[\ga_1(X-\beta)_+-\ga_2(\be-X)_+]\big\},$$ then there exists a constant $\beta_0\in\RR$ such that
$$\beta_0+\EE\big[\gamma_1(X-\beta_0)_+-\gamma_2(\beta_0-X)_+\big]=\EE(X)=\beta_0+\EE\big[(X-\beta_0)_+-(\beta_0-X)_+\big],$$ that is,
$$(\gamma_1-1)\EE[(X-\beta_0)_+]+(1-\gamma_2)\EE[(\beta_0-X)_+]=0.$$ Since $0\le\ga_2<1<\ga_1$, we can get $\EE[(X-\beta_0)_+]=\EE[(\beta_0-X)_+]=0$, and therefore, $X=\beta_0$ almost surely. Hence the OCE measure is averse.

Finally, setting $\ga_1=(1-\al)^{-1}$ and $\ga_2=0$ in \reff{001}, we obtain CVaR. Thus,  CVaR is averse.\qed

\bigskip

On the contrary, we next show that  the risk measure from dividing the future is not averse.
\begin{prop}\label{p:dividing_averse}
The risk  measure defined in (\ref{dividing_formula}) is not averse if  $r\geq2$.
\end{prop}
\pf If $\PP_0(\Omega_k)\neq\lambda_k$ for some $k=1,2,\cdots,r$, then by (\ref{dividing_envelope}), $1\not\in\QQ$. Thus, by Proposition \ref{p:averse_ness}, $\RRR$ is not averse.

If $\PP_0(\Omega_k)=\lambda_k$ for all $k=1,2,\cdots,r$, then set $X=\sum\limits_{k=1}^rk\textbf{1}_{\Omega_k}$. Obviously $X$ is nonconstant. Since
$$\RRR(X)=\sum\limits_{k=1}^r\lambda_k\cdot k=\sum\limits_{k=1}^rk\PP_0(\Omega_k)=\EE(X),$$ which implies that $\RRR$ is not averse.\qed

Although the risk measure from subdividing the future is not averse, this risk measure can be used in composition with other averse measures (say, CVaR) to create new risk measures that make practical sense. We leave this topic for future research.

\section{Coherent risk measures on subspaces: Risk envelopes and uncertainty sets}

Recently, coherent risk measures have been studied in the literature of robust optimization. For instance,  several coherent risk measures were constructed by using the so-called uncertainty sets in Natarajan, Pachamanova, and Sim (2009), while Bertsimas and Brown (2009)  examined the question from a different perspective: If risk preferences are specified by a coherent risk measure, how would the uncertainty set be constructed? In general, from the viewpoint of robust optimization, a risk measure is applied to a random variable of a special structure (say, a linear combination of basic random variables) and is defined by uncertainty sets without involving the exact details of the probability structure of the random variables. In particular, the mean-standard deviation measure, the discrete CVaR, and the distortion risk measure are defined through cone-representable uncertainty sets. If the same risk measure can be constructed by both risk envelope and uncertainty set, then there must be certain relation between the two subjects. It is therefore of interest to explore the connection between risk envelopes and uncertainty sets.  This would help  to have a deeper understanding on robust optimization.

Let us consider a rather general case in robust optimization, where all uncertain data are linear functions of a finite number of random variables, $X_1,...,X_n,$ where $X_i\in\LL^2(\Omega,\Sigma,\PP_0)$ for $1\leq i\leq n$.
Denote
$$\VV:=\left\{X= \rv:~ a_1,\cdots,a_n\in\RR\right\}.$$ Then $\VV$ is the subspace generated by $X_1,...,X_n$.  Let $\RRR$ be a coherent risk measure on $\LL^2(\Omega,\Sigma,\PP_0).$ We define a risk envelope by
\begin{eqnarray}\label{e:envelope}
\QQ_\VV:=&&\left\{Q\in\PPP:~\EE(XQ)\leq\RRR(X)~\h{for all}~X\in\VV\right\}.
\end{eqnarray}
It is easy to check that $\QQ_\VV\su \PPP$ and is nonempty, convex and closed, so it is a risk envelope with an induced risk measure

\begin{equation}\label{e:dual}
\RRR_\VV\left(X\right)=\sup\limits_{Q\in\QQ_\VV}\EE(XQ).
\end{equation}
 Note that the risk envelope $\QQ_\VV$, together with $\RRR_\VV$, relies on the choice of the subspace $\VV$ as well as the original risk measure $\RRR$. Since $\VV$ and $\RRR$ are fixed in the analysis below,  for  notational convenience, we henceforth use $\QQb$ and $\RRb$ for $\QQ_\VV$ and $\RRR_\VV$, respectively. We will also call $\RRb$ the risk measure on $\VV$ to specify its dependence on $\VV$ and $\RRR$.

 We next show that the uncertainty set used in robust optimization for constructing a coherent risk measure on $\VV$ is the (weak) closure  of ``expected image'' of the risk envelope. We need introduce some notations. For any risk envelope $\QQ$, we denote
\begin{equation}\label{e:UQ}
\UU_\QQ:=\cl\left\{\begin{pmatrix}\EE(X_1Q)\cr \vdots\cr \EE(X_nQ) \end{pmatrix}:~Q\in\QQ\right\}.
\end{equation}
In particular, we denote
\[
\UU_\PPP:=\cl\left\{\begin{pmatrix}\EE(X_1Q)\cr \vdots\cr \EE(X_nQ) \end{pmatrix}:~Q\in\PPP\right\}.
\] Then $\UU_\QQ$ is a nonempty and convex subset of $\UU_\PPP$. Given a nonempty,  convex and closed uncertainty set $\UU\subseteq\UU_\PPP$, let

\begin{equation}\label{e:QU}
\QQ_\UU:=\cl\left\{Q\in\PPP:~\begin{pmatrix}\EE(X_1Q)\cr \vdots\cr \EE(X_nQ) \end{pmatrix}\in\UU\right\}.
\end{equation}
Then $\QQ_\UU$ is a nonempty, closed and convex subset of $\PPP$.  The following lemma is basic.
\begin{lemma}\label{property}
The following relations hold:

(1)~$\QQ_{\UU_\PPP}=\PPP$;

(2)~$\UU_{\QQ_\UU}=\UU$;

(3)~$\QQ\subseteq\QQ_{\UU_\QQ}$;

(4)~If $\QQ_1\subseteq\QQ_2$, then $\UU_{\QQ_1}\subseteq\UU_{\QQ_2}$;

(5)~$\UU_1\subseteq\UU_2$ if and only if $\QQ_{\UU_1}\subseteq\QQ_{\UU_2}$.
\end{lemma}
\pf

 (1) Trivial.

(2) On one hand, we have $$\UU_{\QQ_\UU}=\cl\left\{[\EE(X_1Q),..., \EE(X_nQ)]':~Q\in\QQ_\UU\right\}\subseteq\UU,$$ where the  apostrophe stands for the transpose. On the other hand, for any $(z_1,..., z_n)'\in\UU\subseteq\UU_\PPP$, there exists $Q\in\PPP$ such that $z_i=\EE(X_iQ)$ for any $1\leq i\leq n$. Since $[\EE(X_1Q),..., \EE(X_nQ)]'\in\UU$, by definition we have $Q\in\QQ_\UU$. Therefore,
$$(z_1,..., z_n)'=[\EE(X_1Q),..., \EE(X_nQ)]'\in\UU_{\QQ_\UU}.$$ Hence $\UU\subseteq\UU_{\QQ_\UU}$, and then $\UU_{\QQ_\UU}=\UU$.

(3) For any $Q\in\QQ$, we have $[\EE(X_1Q),..., \EE(X_nQ)]' \in\UU_\QQ$. Then by definition, $Q\in\QQ_{\UU_\QQ}$. Therefore, $\QQ\subseteq\QQ_{\UU_\QQ}$.

(4) Trivial.

(5) The ``only if'' part is trivial. For the ``if'' part, by (4) and (2), $\QQ_{\UU_1}\subseteq\QQ_{\UU_2}$ implies $\UU_{\QQ_{\UU_1}}\subseteq\UU_{\QQ_{\UU_2}}$, that is, $\UU_1\subseteq\UU_2$.\qed

\n {\b Remark.} The converse of (3) may not be true. For example, if $\QQ$ is a singleton $\{1\}$, then $\UU_\QQ=[\EE(X_1),..., \EE(X_n)]' $. Here $\QQ_{\UU_\QQ}$ contains all $Q\in\PPP$ such that $[\EE(X_1Q),..., \EE(X_nQ) ]'=[\EE(X_1),..., \EE(X_n) ]'$, which may not necessarily be constant variable $1$.
\bigskip

We can use the uncertainty sets to define coherent risk measures. For uncertainty set $\UU$, the mapping $$ \sum\limits_{i=1}^na_iX_i\longmapsto\sup\limits_{(z_1,\cdots,z_n)^\prime\in\UU}\left( \sum\limits_{i=1}^na_iz_i\right)$$ defines a risk measure on the subspace $\VV$, which is called the risk measure on $\VV$ with uncertainty set $\UU$.

The next two propositions describe some relationships between risk envelopes and uncertainty sets. A common criticism to robust optimization is the arbitrariness of the uncertainty set and its lack of theoretical foundation. Our result here may shed some light on the rationale of uncertainty set and build up a proper theoretical foundation of it. \reft{uncertainty} below serves  for the same purpose.

\begin{prop}\label{rel1}
$\RRb$ is a coherent risk measure on $\VV$ with risk envelope $\QQb$ if and only if it is a coherent risk measure on $\VV$ with uncertainty set $\UU_{\QQb}$.
\end{prop}
\pf By direct calculation, we can get
\beaa
\sup\limits_{Q\in\QQb}\EE\left[\left( \rv\right)Q\right]&=\sup\limits_{Q\in\QQb}\left( \sum\limits_{i=1}^na_i\EE(X_iQ)\right)\\
&=\sup\limits_{(z_1,\cdots,z_n)^T\in\UU_{\QQb}}\left( \rob\right)
\eeaa
for any $ \rv\in\VV$.\qed

\begin{prop}\label{rel2}
For any uncertainty set $\UU\subseteq\UU_\PPP$, $\RRb$ is a coherent risk measure on $\VV$ with uncertainty set $\UU$ if and only if it is a coherent risk measure on $\VV$ with risk envelope $\QQ_\UU$.
\end{prop}
\pf By Proposition \ref{rel1}, $\RRb$ is a coherent risk measure on $\VV$ with risk envelope $\QQ_\UU$ if and only if it is a coherent risk measure on $\VV$ with uncertainty set $\UU_{\QQ_\UU}$. Then by Lemma \ref{property} (2), $\UU_{\QQ_\UU}=\UU$, so the proposition is proved.\qed

 The following is a main theorem in Natarajan $et \ al.$ (2009), where the authors discussed how to construct coherent risk measures in general. However, since uncertainty sets are constructed independent of  probability distributions, it is not completely clear how the uncertainty sets are related to the random variables appeared in the problem. We now present a new proof of the theorem, which discloses the connection between the uncertainty set and the risk measure on $\VV.$
\begin{thm}\label{uncertainty}
$\RRb$ is a coherent risk measure on $\VV$ if and only if there exists a nonempty and convex subset $\UU\subseteq\UU_\PPP$ such that
\begin{equation}\label{e:uncertainty}
\RRb\left( \rv\right)=\sup\limits_{{z}=(z_1,\cdots,z_n)'\in\UU}\left( \rob\right)
\end{equation}
for any $a_1,\cdots,a_n\in\RR$. We call $\UU$ the ``uncertainty set'' of the risk measure $\RRb$ on $\VV.$ It can be written explicitly as $$\UU=\left\{{z}\in\UU_\PPP:~\max\limits_{a_1,\cdots,a_n\in\RR}\left\{\rob:~\RRR\left(\rv\right)\leq1\right\}\leq1\right\},$$where $\RRR$ is the original risk measure that induces $\RRb.$
\end{thm}

\pf Formula \reff{e:uncertainty} follows from Propositions 5.1 and 5.2. 
Next,  by Proposition \ref{rel1},  $\RRb$ is a coherent risk measure on $\VV$ with risk envelope
$$\QQb=\left\{Q\in\PPP:~\EE\left[\left( \rv\right)Q\right]\leq\RRR\left( \rv\right)~\h{for all}~ a_1,\cdots,a_n\in\RR\right\}$$ if and only if it is a coherent risk measure on $\VV$ with uncertainty set
$$\UU_{\QQb}=\left\{\begin{pmatrix}\EE(X_1Q)\cr \vdots\cr \EE(X_nQ) \end{pmatrix}:~Q\in\PPP,~\EE\left[\left( \rv\right)Q\right]\leq\RRR\left( \rv\right)~\h{for all}~ a_1,\cdots,a_n\in\RR\right\}.$$ Therefore, to complete the proof of Theorem \ref{uncertainty}, we only need to prove
\bea\label{equation}
&\left\{\begin{pmatrix}\EE(X_1Q)\cr \vdots\cr \EE(X_nQ) \end{pmatrix}:~Q\in\PPP,~\EE\left[\left( \rv\right)Q\right]\leq\RRR\left( \rv\right)~\h{for all}~ a_1,\cdots,a_n\in\RR\right\}\nonumber\\
=&\left\{\begin{pmatrix}z_1\cr \vdots\cr z_n \end{pmatrix}\in\UU_\PPP:~\max\limits_{a_1,\cdots,a_n\in\RR}\left\{\rob:~\RRR\left(\rv\right)\leq1\right\}\leq1\right\}.
\eea
In fact, since $Q\in\PPP\Longleftrightarrow[\EE(X_1Q),...,\EE(X_nQ)]'\in\UU_\PPP$, and for any $Q\in\PPP$,
\beaa
&&\EE\left[\left( \rv\right)Q\right]\leq\RRR\left( \rv\right)~\h{for all}~ a_1,\cdots,a_n\in\RR\\
&\Longleftrightarrow&\sum\limits_{i=1}^na_i\EE(X_iQ)\leq\RRR\left(\rv\right)~\h{for all}~a_1,\cdots,a_n\in\RR\\
&\Longleftrightarrow&\max\left\{\sum\limits_{i=1}^na_i\EE(X_iQ):~a_1,\cdots,a_n\in\RR,~\RRR\left(\rv\right)\leq1\right\}\leq1,
\eeaa
then (\ref{equation}) holds. The proof of Theorem \ref{uncertainty} is completed.\qed

\section{Concluding Remarks}

Artzner, Delbaen, Eber, and Heath (1997, 1999) introduced the fundamental notion of coherent risk measures. Rockafellar, Uryasev, and Zabarankin (2006)  considered a dual representation theorem in $\LL^2$
space. In this paper,  we considered risk measures in $\LL^2$ under set operations and discussed the dual representations and aversity for various popular risk measures. We also studied the relationship between the risk measure defined by risk envelopes and that defined by uncertainty sets in the case for the risk measures on subspaces. These results may provide  certain tools for stochastic optimization with risk measures as well as  improve our understanding on robust optimization.

\end{document}